\documentclass[12pt,english]{smfarts}
\usepackage{smfthm}
 \usepackage{MnSymbol}
\usepackage{endnotes}

\usepackage[ethiop,english]{babel} \newcommand{\ethi}{\selectlanguage{ethiop}}  
\usepackage{hyperref}

\usepackage{epigraph}
\usepackage{upgreek}

\DeclareMathAlphabet\mathbb{U}{msb}{m}{n}

\usepackage{graphicx}

\newcommand{\bi}{\begin{itemize}}
\newcommand{\ei}{\end{itemize}}
\newcommand{\bd}{\begin{description}}
\newcommand{\ed}{\end{description}}

\newcommand{\bee}{\begin{enumerate}}
\newcommand{\eee}{\end{enumerate}}

\def\hra{\hookrightarrow}

\def\lra{\longrightarrow}
\def\lla{\longleftarrow}
\def\ra{\rightarrow}
\def\llrra{\leftrightarrow}

\def\rtt{\,\rightthreetimes\,}

\makeatletter
\newcommand{\xleftrightarrow}[2][]{\ext@arrow 3359\leftrightarrowfill@{#1}{#2}}
\newcommand{\xdashrightarrow}[2][]{\ext@arrow 0359\rightarrowfill@@{#1}{#2}}
\newcommand{\xdashleftarrow}[2][]{\ext@arrow 3095\leftarrowfill@@{#1}{#2}}
\newcommand{\xdashleftrightarrow}[2][]{\ext@arrow 3359\leftrightarrowfill@@{#1}{#2}}
\def\rightarrowfill@@{\arrowfill@@\relax\relbar\rightarrow}
\def\leftarrowfill@@{\arrowfill@@\leftarrow\relbar\relax}
\def\leftrightarrowfill@@{\arrowfill@@\leftarrow\relbar\rightarrow}
\def\arrowfill@@#1#2#3#4{%
  $\m@th\thickmuskip0mu\medmuskip\thickmuskip\thinmuskip\thickmuskip
   \relax#4#1
   \xleaders\hbox{$#4#2$}\hfill
   #3$%
}

\makeatletter
\newcommand{\xRightarrow}[2][]{\ext@arrow 0359\Rightarrowfill@{#1}{#2}}
\newcommand{\xLeftarrow}[2][]{\ext@arrow 0359\Leftarrowfill@{#1}{#2}}
\makeatother

\makeatletter
\newcommand*{\doublerightarrow}[2]{\mathrel{
  \settowidth{\@tempdima}{$\scriptstyle#1$}
  \settowidth{\@tempdimb}{$\scriptstyle#2$}
  \ifdim\@tempdimb>\@tempdima \@tempdima=\@tempdimb\fi
  \mathop{\vcenter{
    \offinterlineskip\ialign{\hbox to\dimexpr\@tempdima+1em{##}\cr
    \rightarrowfill\cr\noalign{\kern-.3ex}
    \rightarrowfill\cr}}}\limits^{\!#1}_{\!#2}}}
\newcommand*{\triplerightarrow}[1]{\mathrel{
  \settowidth{\@tempdima}{$\scriptstyle#1$}
  \mathop{\vcenter{
    \offinterlineskip\ialign{\hbox to\dimexpr\@tempdima+1em{##}\cr
\rightarrowfill\cr\noalign{\kern-.3ex}
    \rightarrowfill\cr\noalign{\kern-.3ex}
    \rightarrowfill\cr}}}\limits^{\!#1}}}

\newcommand*{\XtoXX}[2]{\mathrel{
  \settowidth{\@tempdima}{$\scriptstyle#1$}
  \mathop{\vcenter{
    \offinterlineskip\ialign{\hbox to\dimexpr\@tempdima+1em{##}\cr
\leftarrowfill\cr\noalign{\kern-.3ex}
    \rightarrowfill\cr\noalign{\kern-.3ex}
    \leftarrowfill\cr}}}\limits^{\!#1}_{\!#2}}}

\newcommand*{\XXtoXXX}[1]{\mathrel{
  \settowidth{\@tempdima}{$\scriptstyle#1$}
  \mathop{\vcenter{
    \offinterlineskip\ialign{\hbox to\dimexpr\@tempdima+1em{##}\cr
\leftarrowfill\cr\noalign{\kern-.1ex}
    \rightarrowfill\cr\noalign{\kern-.3ex}
\leftarrowfill\cr\noalign{\kern-.3ex}
    \rightarrowfill\cr\noalign{\kern-.3ex}
    \leftarrowfill\cr}}}\limits^{\!#1}}}

\newcommand*{\XXXtoXXXX}[1]{\mathrel{
  \settowidth{\@tempdima}{$\scriptstyle#1$}
  \mathop{\vcenter{
    \offinterlineskip\ialign{\hbox to\dimexpr\@tempdima+1em{##}\cr
\leftarrowfill\cr\noalign{\kern-.3ex}
    \rightarrowfill\cr\noalign{\kern-.3ex}
\leftarrowfill\cr\noalign{\kern-.3ex}
    \rightarrowfill\cr\noalign{\kern-.3ex}
\leftarrowfill\cr\noalign{\kern-.3ex}
    \rightarrowfill\cr\noalign{\kern-.3ex}
    \leftarrowfill\cr}}}\limits^{\!#1}}}

\makeatother

\makeatother

\makeatother
\def\xra{\xrightarrow}

\def\ZZ{\Bbb Z}

\def\xra{\xrightarrow}

 \def\inv{^{-1}}

\def\cof{\text{cof}}
\def\lr{\text{lr}}
\def\rl{\text{rl}}
\def\Top{\text{Top}}
\def\preorders{{\text{preorders}}}
\def\Paths{\text{Paths\,}}

\def\antidiscrete{\text{antidiscrete}}
\def\lrl{\text{l}}
\def\rlr{\text{r}}

\def\rrt#1#2#3#4#5#6{\xymatrix{ {#1} \ar[r]^{} \ar@{->}[d]_{#2} & {#4} \ar[d]^{#5} \\ {#3}  \ar[r] \ar@{-->}[ur]^{}& {#6} }}

\def\NN{\Bbb N}
\def\RR{\Bbb R}

\def\card{\,{\mathrm{card}\,}}

\def\AE{\forall\exists}

\def\Dop{\Delta^{\mathrm{op}}}

\def\Sets{\mathrm{Sets}}
\def\sSets{\mathrm{sSets}}

\def\Topp{\mathrm{Top}}

\def\LL{\mathcal L}
\def\NN{\mathbb N}

\def\id{{\text{id}}}
\def\dist{\text{dist}}
\def\diag{\text{diag}}

\def\Filt{{\ethi\ethmath{wA}}}
\def\sFilt{{{\ethi\ethmath{\raisebox{-2.39pt}{nI}\raisebox{2.39pt}\,\!\!wA}}}}

\def\sFilth{\mathit{sFilt}}

\def\Filth{{\ethi\ethmath{wA}}}
\def\sFilth{{{\ethi\ethmath{\raisebox{-2.39pt}{nI}\raisebox{2.39pt}\,\!\!wA}}}}

\def\sFilthLocal{{{\ethi\ethmath{\raisebox{-2.39pt}{nI}\raisebox{2.39pt}\,\!\!wE}}}}

\def\Ob{\text{Ob\,}}

\def\sSet{\text{sSet}}
\def\diag{\text{diag}}
\def\cart{\text{cart}}

\def\const{\text{const}}

\def\FFilt{{\ethi\ethmath{wE}}}
\def\sFFilt{{\ethi\ethmath{\raisebox{-2.39pt}{nE}\raisebox{2.39pt}\,\!\!wE}}}

\def\Stone{{\ethi{\ethmath{cI}}}}

\def\ttt{{\ethi{\ethmath{pa}}}}

\def\mU{{\ethi{\ethmath{mi}}}}

\def\mmU{{\ethi{\ethmath{mA}}}}

\def\Arch{\text{Arch}}

\def\hom#1#2{\left\{#1 \xrightarrow [\text{}]{} #2\right\}}
\def\homm#1#2#3{\left\{{#2} \xrightarrow [\text{}{#1}]{} {#3}\right\}}

\def\hommS#1#2#3{\left\{{#2} \xrightarrow [\text{}{#1}]{sSets} {#3}\right\}}

\def\Hom#1#2{\left\{#1 \xRightarrow [\text{}]{} #2\right\}}
\def\Homm#1#2#3{\left\{{#2} \xRightarrow [\text{}{#1}]{} {#3}\right\}}
\def\hHomm#1#2#3{\left\{{#2} \xRightarrow [\text{}{#1}]{\text{hom}} {#3}\right\}}

\def\HommS#1#2#3{\left\{{#2} \xRightarrow [\text{}{#1}]{sSets} {#3}\right\}}

\def\homSets#1#2{\text{Hom}\left(#1 ,  #2\right)}
\def\hommSets#1#2#3{\text{Hom}_{#1}\left({#2} , {#3}\right)}
\def\HomSets#1#2{\text{\underline{Hom}}\left(#1 ,  #2\right)}
\def\HommSets#1#2#3{\text{\underline{Hom}}_{#1}\left({#2} , {#3}\right)}
\title[Geometric realisation as a space of semi-continuous paths]
{A geometric realisation of\\ geometric realisation\\ as the Skorokhod semi-continuous path space endofunctor
}
\author{Misha Gavrilovich, Konstantin Pimenov}
\address{a draft. $\bullet$ m.gavrilovich 
National Research University Higher School of Economics, Saint-Petersburg;
Institute for Regional Economics Studies of the Russian Academy of Sciences (IRES RAS)                                                                       
38 Serpuhovskaya st., Saint-Petersburg.} 
\address{konstantin pimenov St.Petersburg State University}

\email{mi\!\!\!ishap\!\!\!p@sd\!\!\!df.org}
\urladdr{http://mishap.sdf.org/SkorokhodGeometricRealisation.pdf
{http://mishap.sdf.org/SkorokhodGeometricRealisationHomSets.pdf}
{http://mishap.sdf.org/SkorokhodGeometricRealisation.tex}
{\ \ \ \ Discussion group \href{https://t.me/joinchat/GVRrKxbSO8EWehZYReTKeQ}{here}.}
}
\tiny\thanks{
\tiny{\tiny 
You have a soul of your own, and the privilege of
free-will, whoever you be, as well as the proudest he that struts in a gaudy
outside; you are a king by your own fireside, as much as any monarch on his
throne; you have liberty and property, which set you above favour or affection,
and may therefore freely like or dislike this history, according to your
humour.} 
}

\begin{document}\selectlanguage{english}\catcode`\_=8\catcode`\^=7 \catcode`\_=8

\begin{abstract} We interpret a construction of geometric realisation by [Besser], [Grayson], and [Drinfeld] 
as constructing a space of maps from the interval to a simplicial set, in a certain formal sense, 
reminiscent of the Skorokhod space of semi-continuous functions; 
in particular, we show the geometric realisation functor factors through an endofunctor of a certain category. 
Our interpretation clarifies
the explanation of [Drinfeld] ``why geometric realization commutes with Cartesian
products and why the geometric realization of a simplicial set [...] 
is equipped with an action of the group of orientation preserving
homeomorphisms of the segment $[0, 1]$''. 
\end{abstract}
%
%
%
%
%
%
%
\maketitle
{\small\tiny 
\setcounter{tocdepth}{2}
\tableofcontents
}

\section
{Introduction}

%
%

We show that the geometric realisation functor $|\cdot|:\sSets\lra\Topp$ factors via an endofunctor of 
a certain category\footnote{We suggest to pronounce $\sFilth$ as $sF$
as it is visually similar to 
$s\Phi$ standing for ``{\em s}implicial $\phi$ilters'', 
even though it is unrelated to the actual pronunciation of these symbols coming from the Amharic script.
} 
 $\sFilth$ extending both the categories of simplicial sets and topological spaces (as full subcategories); 
this endofunctor can be thought of,  in a certain formal sense,  as 
the space of semi-continuous maps from the interval to a simplicial set, 
reminiscent of the Skorokhod space used in probability theory.
The category  $\sFilth$ can be thought of as the category of simplicial sets
equipped with additional structure of topological flavour, 
namely the notion of a simplex ``being small enough''.
Formally, the endofunctor is the following inner hom
of 
simplicial sets equipped with a certain additional structure (see \S\ref{sketch}): 
$$ \Homm{sSets}{\homm{\text{preorders}}{-}{[0,1]_\leq}}{ {X_\bullet} } = 
 \HommSets{sSets}{\hommSets{\text{preorders}}{-}{[0,1]_\leq}}{ {X_\bullet} }
$$

We show this by rewriting in terms of $\sFilth$ a construction of [Besser], [Grayson], and [Drinfeld],
and thus the words of [Drinfeld] apply to our paper as well (with the
exception that we do not discuss cyclic sets):
\begin{quote}
We explain why geometric realization commutes with Cartesian 
products and why the geometric realization of a simplicial set [..] 
is equipped with an action of the group of orientation preserving 
homeomorphisms of the segment $[0, 1]$.

....

In this note there are no theorems, its only goal is to clarify the notion of
geometric realization for simplicial sets [...] 
We reformulate the definitions so that the following facts become obvious:
\bi\item[(i)] geometric realization commutes with finite projective limits (e.g., with
Cartesian products);
\item[(ii)] the geometric realization of a simplicial set [...] 
is
equipped with an action of the group of orientation preserving homeomorphisms of the segment I := [0, 1] [...].
\ei
In the traditional approach [...] 
these statements are
theorems, and understanding their proofs requires some efforts.
\end{quote}


We aim to be brief, self-contained, and include only what is needed formally; 
a reader so inclined may look up a verbose discussion of constructions and intuitions in preliminary notes 
\href{http://mishap.sdf.org/6a6ywke/6a6ywke.pdf}{[6]}. Our definitions are
complete and self-contained; proofs are routine verifications sometimes
omitted.

The category $\sFilth$ is the category of simplicial objects in the category of filters
in the sense of Bourbaki (defined in \S\ref{filter:def}).
It is rather large: it contains, as full subcategories, in several ways,
the categories of topological and metric (uniform) spaces; the category of simplicial sets, 
and the category of filters. As the category contains both filters and spaces,
the notion of {\em limit} is naturally expressible in it, in fact 
in terms of the lifting properties and 
the {\em ``shift''} endofunctor $[+1]:\Delta\lra\Delta,n\longmapsto n+1$ of the category of finite linear orders. 
The same endofunctor can be used to reformulate the standard definition 
of  {\em locally trivial fibre bundle} of topological spaces, 
as a morphism  becoming a trivial fibre bundle (i.e.~a direct product)
after a base-change. 
In both examples,  the ``shift'' endofunctor  $[+1]:\Delta\lra\Delta$ is 
used to talk about local properties. One may associate  $\sFilth$-objects
to  metric spaces  in such a way that  quasi-isometries
become $\sFilth$-isomorphisms.

We sketch some of these preliminary definitions and constructions in Appendix~\S\ref{remarks}, 
as any discussion of them is outside the scope of this short note. We want to make an explicit warning
that Appendix~\S\ref{remarks} is not ready for publications but we release 
this unfinished work as it is, 
in the small hope that some readers may find it useful.

%
%
%
%
%
%
%


Let us now explain how our construction may arise from an attempt to reformulate a construction of [Grayson] and [Drinfeld] in terms of 
the category $\sFilth$ of simplicial filters. Although formally unnecessary, the reader may want to read first 
the definition of $\sFilth$ and our notation in \S\ref{filter:def},

\subsubsection*{Notation for homs}
To aid readability of nested homs\footnote{\href{http://mishap.sdf.org/Skorokhod_Geometric_Realisation_HomSets.pdf}{Here} the reader may find this paper in the usual notation.
}, we use Haskell-like notation $\hom X Y:=\homSets X Y$ and $\homm {C} X Y:= \hommSets C X Y$  to denote the set of morphisms from an object $X$ to an object $Y$ in a category $C$,
and  $\Hom X Y:=\HomSets X Y$ and $\Homm {C} X Y:= \HommSets C X Y$ to denote the inner hom. 
Thus, for example, 
$$ \homm{sSets}{\homm{\text{preorders}}{-}{[0,1]_\leq}}{ {X_\bullet} } = 
 \hommSets{sSets}{\hommSets{\text{preorders}}{-}{[0,1]_\leq}}{ {X_\bullet} }
$$
$$ \Homm{sSets}{\homm{\text{preorders}}{-}{[0,1]_\leq}}{ {X_\bullet} } = 
 \HommSets{sSets}{\hommSets{\text{preorders}}{-}{[0,1]_\leq}}{ {X_\bullet} }
$$

\subsection*{Our construction: a sketch\label{sketch}} 
First we observe that the set of points of geometric realisation of a simplicial set $X_\bullet$, as defined by  [Drinfeld,\S1, Eq.(1.1)] 
(see \S\ref{drinfeld}), is
$$ \homm{sSets}{\homm{\text{preorders}}{-}{[0,1]_\leq}}{ {X_\bullet} } $$
thereby making it apparent that the set of points of the geometric realization so defined
commutes with projective limits, 
and is
equipped with an action of the group of orientation preserving homeomorphisms of the segment $ [0, 1]_\leq$. 
For $X_\bullet:=\Delta_N=\homm{\text{preorders}}{-}{(N+1)_\leq}$ the $N$-dimensional simplex, by Yoneda lemma 
it is the set of non-decreasing functions $[0,1]\lra (N+1)_\leq$:
$$ \homm{sSets}{\homm{\text{preorders}}{-}{[0,1]_\leq}}{ {\homm{\text{preorders}}{-}{(N+1)_\leq}} } =
\homm{\text{preorders}}{[0,1]_\leq}{ (N+1)_\leq}
$$
As noted by by  [Grayson, 2.4.1-2] (see \S\ref{grayson:remark}), it becomes the standard simplex 
$$|\Delta_N|:=\{ (t_0, ..., t_N)\in \RR^{N+1}\,:\,0\leq t_0\leq ..\leq t_N\leq 1\}$$ in $\RR^{N+1}$ 
when equipped with the Skorokhod metric originating in probability theory  [Kolmogorov, \S2, Def.1] 
$$\dist(f,g):= \inf\{\varepsilon: \forall x \exists y \,(\, |x-y|<\varepsilon \,\&\, |f(x)-g(y)|<\varepsilon\,)\, \}$$
An ordered pair $x\leq y$, $x,y\in[0,1]$ is a $1$-simplex of $\homm{\text{preorders}}{-}{[0,1]_\leq}$,
and an ordered pair of functions $f\leq g$, $\forall x\in [0,1]\,\, ( f(x)\leq g(x))$, is a $1$-simplex of 
the inner hom of simplicial sets
$$ \Homm{sSets}{\homm{\text{preorders}}{-}{[0,1]_\leq}}{ {\homm{\text{preorders}}{-}{(N+1)_\leq}} } $$
Thinking intuitively of these ordered pairs/simplicies of nearby points/functions as 
small perturbations/homotopies $x\rightsquigarrow y$ and $f\rightsquigarrow g$,
or perhaps rather  $(x,x)\rightsquigarrow(x,y)$, $(x,x,x)\rightsquigarrow (x,x,y) \rightsquigarrow (x,y,y)$, and 
$(f,f)\rightsquigarrow(f,g)$, $(f,f,f)\rightsquigarrow (f,f,g) \rightsquigarrow (f,g,g)$, 
and that being Skorokhod-nearby means that 
some small enough perturbation/homotopy of the argument results in a small perturbation/homotopy of the value
(see~\S\ref{skoro:filter}) 
in \S\ref{skoro:space} we define a generalisation of the
Skorokhod metric which equips the inner hom with an extra structure, a notion of smallness, 
thereby turning it into an object of $\sFilth$ which remembers the topology of the geometric realisation, 
and in fact such that a certain forgetful functor $\sFilth\lra \Topp$ (see \S\ref{sFtoTop}) turns it into 
the topological space $|X_\bullet|$, under some finiteness assumptions on the simplicial set $X_\bullet$ 
(see \S\ref{sFDelta2Top}).

\section{Main construction}

\subsection{\label{filter:def}The category of simplicial filters}
A {\em filter} on a set $X$ is a set $\mathfrak F $ of subsets of $X$ which has the following properties :
\bi
\item[$\text{(F}_{\text{I}}\text)$] Every subset of $X$ which contains a set of $\mathfrak F$	belongs to $\mathfrak F$.

\item[$\text{(F}_{\text{II}}\text)$] Every intersection of two subsets of $\mathfrak F$ belongs to $\mathfrak F$.
\item[$\text{(F}'_{\text{II}}\text)$]  $X\in\mathfrak F$.

\ei
Subsets in $\mathfrak F$ are called {\em neighbourhoods} or {\em $\mathfrak F$-big}.
By abuse of language, often by a filter we mean a set together with a filter on it.

A {\em morphism of filters} is a mapping of underlying sets 
such that the preimage of a neighbourhood is necessarily a neighbourhood;
we call such maps of filters {\em continuous}.  

Unlike the definition of filter in \href{http://mishap.sdf.org/tmp/Bourbaki_General_Topology.djvu#page=63}{[Bourbaki, I\S6.1, Def.I]},
we  do not require 
that $\emptyset\not\in\frak F$, in particular it is possible that $X=\emptyset$. 
We do so in order for the category  of filters to have  limits. 

Let $\Filt$ denote the category of filters. 

Let  $\sFilt=Func(\Dop, \Filt)$ be 
the category of functors
from 
$\Dop$, the category opposite to the category $\Delta$ of finite linear orders,
to the category $\Filt$ of filters. One may refer to an object of $\sFilth$
as either {\em a simplicial filter} 
or {\em a situs},
for lack of a short 
name. 

Though we will not be using it, let $\sFilthLocal:=\sFilth/\!\!\approx$ be $\sFilth$ where we consider two functors equal as morphisms 
iff they coincide on a big subset (equivalently, in another terminology, locally coincide): that is,  for $f,g:X_\bullet\lra Y_\bullet$,
\bi\item $f=_\sFilthLocal g$ iff there is $N>0$ such that for each $n>N$ there is a neighbourhood in the set of $n$-simplicies of $X$
such that for each $s\in \varepsilon$ $f(s)=g(s)$.
\ei

\subsubsection*{Our notation}
For a natural number $k\in \NN$, by $k_\leq $ we denote the linearly ordered set $1<2<..<k$ with $k$ elements 
viewed as an object of $\Dop$; a morphism $\theta:k_\leq \lra l_\leq $ in $\Delta$ 
we denote by $[i_1\leq ... \leq i_k]$ where $1\leq  i_1=\theta(1)\leq ... i_k=\theta(k)\leq l $, and 
 $[i..j]$ is short for $[i\leq i+1\leq ..\leq j]$; hence, 
for an object $X_\bullet$ of $\sFilth$ and an $l$-simplex $x\in X_l$,  
its face and degeneracy maps are denoted by $x[i_1\leq ... \leq i_k]=\theta^*(x)\in X_k$,
and $x[i..j]\in X_k$ is short for $x[i\leq i+1\leq ..\leq j]$. 
Sometimes we write $X_\bullet:\sFilth$
to indicate that $X_\bullet$ is an object of $\sFilth$. 
By $\hom X Y:=\homSets X Y$ or $\homm {C} X Y:=\hommSets C X Y$  we denote the set of morphisms from $X:C$ to $Y:C$ in a category $C$; 
hence both  $\hom - {N_\leq}$ and $\homm {\text{preorders}}  - {N_\leq}$ denote the $(N-1)$-dimensional simplex (as a simplicial set). 
By $\Hom X Y:=\HomSets X Y$ or $\Homm {C} X Y:=\HommSets C X Y$  we denote the inner hom from $X:C$ to $Y:C$ in a category $C$
whenever it is defined.

For a simplicial object $X_\bullet$ of a category and $n>0$, by $X_\bullet(n_\leq)$, $X(n_\leq)$, and $X_{n-1}$,  
we denote the object of $n-1$-simplicies; thus for example $X_0=X(1_\leq)=X_\bullet(1_\leq)$. 

Thus $\Homm {\sSets} X Y$ denotes the inner hom in the category of simplicial sets, and thereby
 $$\Homm {\sSets} X Y(1_\leq)={\Homm {\sSets} X Y}_{0}=\homm{\sSets} X Y=\hommSets{\sSets} X Y $$
denotes the set of morphism from $X$ to $Y$. 

We put in quotation marks words intended to aid intuition but formally unnecessary; thus formally 
 $x\in \delta$ and  ``$\delta$-small'' $x\in \delta$ mean the same. 



\subsection{The Skorokhod filter on set of morphisms}\label{skoro:filter}\label{skorokhod-filter-def}\label{Skorokhod-mapping-spaces-sample}


In probability theory, it is natural to consider that a small distortion of either times   or values
of a stochastic process does not change it much; this is captured by the notions of
 Levi-Prokhorov or Skorokhod metric or convergence [Kolmogorov, \S2, Def.1 of $\varepsilon$-equivalence].

Intuitively, a function is {\em Skorokhod small} if, 
 whenever the argument
is sufficiently small, 
its value  can be made as small as we please
by making {\em some} perturbation (``homotopy'') of the argument 
as small as we please.

Let $X_\bullet,Y_\bullet:\sFilth$ be objects of $\sFilth$. We shall now define 
{the Skorokhod filter on the set of morphisms (in sSets)
 $\homm\sSets  {X_{\bullet}} {Y_{\bullet}}$}
of the underlying simplicial sets 
of $X_\bullet$ and $Y_\bullet$. 

For $N>2n$, neighbourhoods $\delta\subset X_N$ and $\varepsilon\subset Y_n$, 
a {\em $\varepsilon\delta$-Skorokhod neighbourhood in 
 Hom-set $\homm\sSets  {X_{\bullet}} {Y_{\bullet}}$} of the underlying simplicial sets
of ${X_\bullet}$ and ${Y_\bullet}$
 is the subset 
consisting 
of all the functions $f:X_{\bullet}\lra Y_{\bullet}$ with the following property:
\bi\item[] there is a neighbourhood $\delta_0\subset X_n$ such that 
each ``$\delta_0$-small'' $x\in \delta_0$ has a ``$\delta$-small''
``continuation'' $x'\in \delta\subset X_N$, $x=x'[1..n]$ such that its ``tail'' maps into something ``$\varepsilon
$-small'', 
i.e.~$f(x'[N-n+1..N])\in\varepsilon$.
\item [] As a formula, this is 
$$\{ f:{X_\bullet}\lra {Y_\bullet}\,:\, \exists \delta_0 \subset X_n\, \forall x\in\delta_0 \, \exists x'\in \delta
\,\,\, (\,\, x=x'[1..n] \,\,\,\&\,\,\, f(x'[N-n+1..N])\in \varepsilon \,) \}$$
\ei
%

Let $\hommS\Filt {X_\bullet} {Y_\bullet}$ denote the set $\homm\sSets  {X_{\bullet}} {Y_{\bullet}}$
equipped with the {\em Skorokhod} filter generated by all the 
 $\varepsilon\delta$-Skorokhod neighbourhoods
where  $N>2n$, and $\delta\subset X_N$  and $\varepsilon\subset Y_n$
are neighbourhoods.

For each ${X_\bullet}:\sFilth$ this defines a functor $$\hommS\Filt {X_\bullet} - : \sFilth \lra \sFilth$$

For a morphism $f:X_\bullet'\lra X_\bullet''$ in $\sFilth$, the induced map 
$$\hommS\Filt {X_\bullet'} {Y_\bullet}\lla \hommS\Filt {X_\bullet''} {Y_\bullet}$$
is continuous whenever the filters on $X_\bullet'$ are induced from $X_\bullet''$ along $f$, 
i.e.~for each $n\geq 0$  the filter on ${X}'_n$ is the coarsest filter such that the map $f_n:{X}'_n \lra {X}''_n$
is continuous.

The {\em Skorokhod filter on  $\homm\sSets  {X_{\bullet}} {Y_{\bullet}}$} 
is the filter generated by all the Skorokhod $\varepsilon\delta$-neighbourhoods for 
$N\geq 2n>0$ (sic!), neighbourhoods $\delta\subset X_N$ and $\varepsilon\subset Y_n$. 
We denote it by 
 $\hommS\Filt {X_\bullet} {Y_\bullet}$.

\subsection{\label{skoro:space}The Skorokhod space of maps}

For $N>2n$, a sequence $\theta=[l_1\leq ... \leq l_n]\in \homm{\Delta}{n_\leq}{M_\leq}$ where $1\leq l_1\leq l_n\leq M$, 
and neighbourhoods $\delta\subset X_N$ and $\varepsilon\subset Y_n$, 
a {\em $\theta\varepsilon\delta$-Skorokhod neighbourhood of 
 Hom-set $\homm\sSets {\hom{-}{M_\leq} \times {X_{\bullet}}} {Y_{\bullet}}$} of the underlying simplicial sets
of ${X_\bullet}$ and ${Y_\bullet}$
 is the subset 
consisting 
of all the functions $f: \hom-{M_\leq} \times X_{\bullet}\lra Y_{\bullet}$ with the following property:
\bi\item[] there is a neighbourhood $\delta_0\subset X_n$ such that 
each ``$\delta_0$-small'' $x\in \delta_0$ has a ``$\delta$-small''
``continuation'' $x'\in \delta\subset X_N$, $x=x'[1..n]$ such that its ``tail'' maps into something ``$\varepsilon
$-small'', 
i.e.~$f(\theta, x'[N-n+1..N])\in\varepsilon$.
\item [] As a formula, this is 
$$\{ f:{X_\bullet}\lra {Y_\bullet}\,:\, \exists \delta_0 \subset X_n\, \forall x\in\delta_0 \, \exists x'\in \delta
\,\,\,\,(\,\, x=x'[1..n] \,\,\,\&\,\,\, f(\theta,x'[N-n+1..N])\in \varepsilon\,\, ) \}$$
\ei

This notion of a neighbourhood turns the inner hom $\Homm\sSets  {X_{\bullet}} {Y_{\bullet}}$
of the underlying simplicial sets into a simplicial filter
which we call {\em the Skorokhod space of maps from $X_\bullet$ to $Y_\bullet$} and denote by 
 $\HommS\Filt {X_\bullet} {Y_\bullet}$.

For each ${X_\bullet}$ this gives a functor $$\HommS\Filt {X_\bullet} - : \sFilth \lra \sFilth$$

For a morphism $f:X_\bullet'\lra X_\bullet''$ in $\sFilth$, the induced map
$$\HommS\Filt {X_\bullet'} {Y_\bullet}\lla \HommS\Filt {X_\bullet''} {Y_\bullet}$$
is continuous whenever the filters on $X_\bullet'$ are induced from $X_\bullet''$ along $f$.

\subsection{The main diagonal and subdivision filter on a simplicial set}

\subsubsection{The filter of the main diagonal} 
%
%
%

Given a simplicial set $X_\bullet:\sSet$, equip $X_0$ with the {\em antidiscrete} filter $\{X_0\}$, 
i.e.~the filter on $X$ consisting of the unique neighbourhood $X_0$ itself;
for each $n>0$ equip $X_n$ with the finest filter such that
the ``main diagonal'' degeneracy map $X_0\lra  X_n$ is continuous; explicitly,
a subset of $X_n,n\geq 0$ is a neighbourhood iff it contains the image of 
the diagonal map $X_0\lra X_n$. Denote by $X_\diag$ the simplicial filter obtained. 
This defines a fully faithful embedding $$\diag:\sSets\lra\sFilth, \,\,\, X\longmapsto X_\diag$$ 
of the category of simplicial sets.

\subsubsection{The subdivision filter on a simplicial set}

%
Let ${X_\bullet}:\sSet$ and let $\epsilon\in X_k$ be a simplex, $m\geq 0$. 
Call a subset $\varepsilon\subset X_n$ is an {\em  $(\epsilon,m)$-neighbourhood}, 
iff
\bi 
\item the face $\epsilon'[t_1\!\leq...\leq\!t_n]$ of simplex $\epsilon'$ is $\varepsilon$-small, i.e.~$\epsilon'[t_1\!\leq...\leq\!t_n]\in \varepsilon$,
 whenever 
\bi\item 
$N\geq 0$ and 
 $0\leq t_1\leq ... \leq t_n\leq t_1+m\leq N$ 
\item $\epsilon'$ is a simplex of dimension $N$ 
\item $\epsilon$ is a face of $\epsilon'$ 
\ei
\ei

The {\em subdivision 
 filter} on $X_n$ is generated by 
  the $(\epsilon,m)$-neighbourhoods where $\epsilon:X_k$ varies through simplices of arbitrary dimension  $k\geq 0$, $m>0$.
{\em The subdivision filter on a simplicial set $X_\bullet$} is the simplicial set $X_\bullet$ together 
with the subdivision filter on $X_n$ for each $n\geq 0$.

Intuitively, a subset is a neighbourhood iff it contains all simplicies which 
``fit in'' ``between'' {\em consecutive} faces of a simplex (of sufficiently high dimension).

Let $X_{\text{subd}}$ denote the simplicial set equipped with the subdivision neighbourhood structure.

\subsection{The unit interval object}\label{unit-interval} 

Let  $[0,1]_\leq$ denote 
the simplicial set 
 co-represented by $[0,1]_\leq$ in the category of preorders  
$$ n_\leq\longmapsto \homm{{\text{preorders}}}{n_\leq}{[0,1]_\leq}$$ 
equipped with the the subdivision filters.
Explicitly, an $n$-simplex is a non-decreasing sequence $t_0\leq ... \leq t_n$,
and the subdivision filter on the set of $n$-simplicies
is generated by 
``$\epsilon$-neighbourhoods of the main diagonal''
$$\{(t_1\leq ...\leq t_n): \dist(t_i,t_{i+1})<\epsilon\text{ for all } 1\leq i\leq i+1\leq n\,\}$$ where $\epsilon>0$.

\subsection{Geometric realisation as a Skorokhod space of paths} 

Now we are ready to consider the Skorokhod space of ``semi-continuous'' paths 
defined in \S\ref{skoro:space}
$$\HommS\Filth{[0,1]_\leq} {X_\diag}$$

We  now proceed to analyse this space, define a forgetful functor to the  category of topological spaces,
and see that it gives the geometric realisation when applied to this space. 

First we analyse this space for the simplex $X_\bullet= \Delta_{N-1}=\hom-{N_\leq}$:
$$\HommS\Filth{[0,1]_\leq} {\hom-{N_\leq} }$$

\subsubsection{\label{grayson:remark}The standard simplex in $\RR^n$ as a space of semi-continuous functions with a Levi-Prokhorov kind of metric} 

As remarked in [Grayson, 2.4.1-2], 
a non-decreasing map $f:[0,1]_\leq \lra N_\leq$ 
 is almost determined (say, up to finitely many points) by 
the sequence $0\leq s_1\leq ... \leq s_{N-1}\leq 1$ of points where it it is discontinuous,
i.e.,~a point of the geometric realisation $|\Delta_{N-1}|$ of the $(N-1)$-simplex $\Delta_{N-1}=\hom-{N_\leq}$.
Conversely, taking connected components of the interval with these points punctured out gives rise to a map to finite linear order. 
 In notation,  
$$|\Delta_{N-1}|=\{(s_1,..,s_{N-1})\in \RR^{N-1}: 0\leq s_1\leq ... \leq s_{N-1}\leq 1\}
 \approx                                                                             
\{ s\!:\![0,1]_\leq \lra N_\leq \} $$ $$                                                                                                                   
 0\leq s_1\leq ... \leq s_{N-1}\leq 1 \ \approx (  [0,s_1)\mapsto 0, ..., [s_{N-2},s_{N-1})\mapsto N-1, [s_{N-1},1]\mapsto N ) 
$$ $$ 
\,\,\,\,\,\,\,\,\,\,\,\,\,\,\,\,\,\,\,\,\,\,\,\,\,\,\,\,\,\,\,\,\,\,\,\,\,\,\,\,  \approx  \,\,\left(\,\, [0,1]\lra \,\pi_0([0,1]\setminus\{s_1,...,s_N\})\hra N \,\,\right)
$$

In terms of functions, the $\sup$-metric on $\Delta_{N-1}$ is reminiscent of Skorokhod 
or Levi-Prokhorov kind of metric well-known in probability theory:
$$\dist_{\sup}( (s_1\leq ... \leq s_{N-1}), (t_1\leq ... \leq t_{N-1})) = \sup_i |s_i-t_i|= $$ $$
=\inf \{ \delta>0: \,\forall \,0\leq t\leq 1 \,\, \exists \, 0\leq t'\leq 1 \text{ such that } |t'-t|  \leq \delta \text{ and } |f(t')-g(t')|<\delta \} $$

\subsubsection{The Skorokhod spaces of paths in a simplex} 
By Yoneda lemma, to give a map  $\hom-{[0,1]_\leq} \lra {\hom{-} {{N_\leq}}}$
of the simplicial sets co-represented by linear orders $[0,1]_\leq$ and $N_\leq$
is to give a non-decreasing map $f:[0,1]_\leq \lra N_\leq$.

To give a map of ssets $\hom{-} {[0,1]_\leq} \times  {\hom{-} {{M_\leq}}} \lra {\hom{-} {{N_\leq}}}$ 
is to give an increasing sequence of functions  $f_1\leq ...\leq f_M:[0,1]_\leq \lra N_\leq$ 
(where $f_i\leq f_j$ means for all $t\in [0,1]$ $\, f_i(t)\leq f_j(t)$):
$$\hom{n_\leq} {[0,1]_\leq} \times  {\hom{n_\leq} {{M_\leq}}} \lra {\hom{n_\leq} {{N_\leq}}}$$
$$(0\leq t_1\leq ... \leq t_n\leq 1)\times (1\leq l_1\leq ... \leq l_n\leq M)
\longmapsto 
(f_{l_1}(t_1) \leq ... \leq f_{l_n}(t_n))
$$

The Skorokhod filter on $\hom{-} {[0,1]_\leq} \times  {\hom{-} {{M_\leq}}} \lra {\hom{-} {{N_\leq}}}$ 
can be explicitly described as follows. 

First note that we may assume  that $\varepsilon\subset \hom{n_\leq } {{N_\leq}} $ is the smallest  possible, 
i.e. the main diagonal $\varepsilon=\{(i\leq i\leq ..\leq i): 1\leq i\leq  N\}$,
and that $\delta=\{(t_1\leq ... \leq t_n: 0\leq t_1\leq t_n\leq t_1+\updelta, t_n\leq 1\}$ for some $\updelta>0$. 
The corresponding Skorokhod neighbourhoods for $N>2n$ and a sequence 
$(1\leq l_1\leq ... \leq l_n\leq M)$ consists of the increasing sequences of functions
$f=(f_1\leq ... \leq f_M):[0,1]_\leq \lra N_\leq$ such that
\bi\item[] there is $\updelta_0>0$  such that 
each ``$\updelta_0$-small'' sequence $0\leq t_1\leq...\leq t_n\leq t_1+\updelta_0$ 
has a ``$\updelta$-small''
``continuation'' $0\leq t_1\leq...\leq t_n\leq t_{n+1}\leq ...\leq t_N\leq t_1+\updelta$, $t_N\leq 1$, 
such that its ``tail'' maps into something ''on the main diagonal'' 
i.e.~$f_{l_1}(t_{N-n+1})= ... =f_{l_N}(t_N)$.
\ei

Note that when $M=1$, we may pick $t_{N-n}= ... =t_N$ all equal, and therefore 
the Skorokhod filter on the $0$-simplicies of the Skorokhod space is antidiscrete. 

For $M=2$, the condition above says that for any $0\leq t\leq 1$  $f_2(t')=f_1(t')$ for some $t\leq t'\leq t+\delta$ and $t'\leq 1$.
More generally, for any $M>0$, the condition above says that for any $0\leq t\leq 1$ 
 $f_1(t')=...=f_M(t')$ for some $t\leq t'\leq t+\delta$ and $t'\leq 1$.      

%

\subsubsection{\label{sFtoTop}Forgetful functor to topological spaces}  

We shall now apply to the Skorokhod space the ``forgetful'' functor 
 $\sFilth\lra \Topp$ defined as follows. 
For a simplicial filter ${X_\bullet}:\sFilth$,
let 
$$X_\text{points} := \bigcap\limits_{\varepsilon\subset X_0\text{ such that } \varepsilon\text{ is a neighbourhood}} \varepsilon$$
be the set of points of the corresponding topological space,
and the topology be the coarsest topology such that 
for each point $x\in X_\text{points}$ and each neighbourhood $\delta\subset X_1=X(2_\leq)$,
the set of points 
$$U_{\delta,x}:=\left\{y: x=\gamma[0]\text{ and }\gamma[1]=\gamma'[1]\text{ and }
\gamma'[0]=y\text{ for some }\gamma,\gamma'\in\delta\,\right\}
$$
is a neighbourhood (``the $\delta$-neighbourhood of $x$''). 

More graphically, 
$$U_{\delta,x}:=\left\{y: x \xra\gamma \cdot \xleftarrow {\gamma'} y \text{ for 
 some }\gamma,\gamma'\in\delta\,\right\}
$$
is a neighbourhood (``the $\delta$-neighbourhood of $x$'').

 
\subsubsection{\label{sFDelta2Top}The topological space associated with the Skorokhod space of paths of a simplex.} 
Let us now compute the 
topological space associated with the Skorokhod space of paths of a standard simplex.
For the Skorokhod space of paths, 
this construction  gives 
$X_\text{points} = \hom-{[0,1]_\leq} \lra {\hom{-} {{N_\leq}}}$
is the set of non-decreasing functions $f:[0,1]\lra N_\leq$. 
A $1$-simplex is an increasing pair of non-decreasing functions $f_1\leq f_2:[0,1]\lra N_\leq$,
and it lies in $\theta\delta$-Skorokhod neighbourhood for each $\theta:2\lra 2 $ iff
for each $t$ there is $t'\leq t+\delta$ such that $f_2(t')=f_1(t)$ (and $t'\leq 1$). 

A subset is a Skorokhod neighbourhood of a non-decreasing function $f:[0,1]\lra N_\leq$
iff there is $\delta>0$ such that it contains  all the 
non-decreasing functions $f':[0,1]\lra N_\leq$
such that 
\bi\item for every $0\leq t\leq 1$ $f'(t')=f(t')$ for some $0\leq t'\leq 1$
such that $t\leq t'\leq t+\delta$. 
\ei

To see this, consider the function $max(f,f')$ 
and note that both $1$-simplicies 
$(f\leq \max(f,f'))$ are well-defined and  
lie in a $\theta\delta$-Skorokhod neighbourhood for each $\theta$. 

Note that the functions corresponding to the same sequence $0\leq s_1\leq ... \leq s_{N-1} \leq 1$
are topologically indistinguishable. 
In terms of the sequences, the definition above is rephrased as follows: 
\bi\item 
 $|s'_i-s_i|\leq \delta$ for all $1\leq i \leq N-1$. 
\ei
and thus we see this is the usual topology of the simplex in $\RR^{N-1}$
defined by the $\sup$-metric, up to indistinguishable points.

This finishes the interpretation of the construction of [Grayson, Remark 2.4.1-2] 
of the geometric realisation in terms of $\sFilth$.

\subsubsection{\label{drinfeld}The Skorokhod space of an arbitrary simplicial set} 

Now we follow the construction of [Drinfeld,\S1]. 

A  finite subset $F\subset [0,1]$ and an $x\in {X_\bullet}(\pi_0([0,1]\setminus F))$ determines a morphism of sSets 
$ \homm{\text{preorders}}{-}{[0,1]_\leq} \lra  {X_\bullet} $ as follows:
$$ \homm{\text{preorders}}{n_\leq}{[0,1]_\leq} \lra  {X_\bullet}(n_\leq) $$
$$\overrightarrow t \longmapsto x\left[\xra{n_\leq\xra{\overrightarrow t} [0,1]\lra \pi_0([0,1]\setminus F)}\right]\in {X_\bullet}(n_\leq)
$$
where $[0,1]\lra \pi_0([0,1]\setminus F)$ is the obvious map contracting the connected components (we need to make
a convention where to send points of $F$).

According to [Drinfeld, (1.1)], the geometric realisation can be defined as 
$$                                                                                                                                                             
|{X_\bullet}|:=\varinjlim\limits_{{ F\subset [0,1] \text{ finite}}} {X_\bullet}(\pi_0([0,1]\setminus F))$$
where $F$ runs through the set of all finite subsets of $[0,1]$ and $\pi_0([0,1]I \ F )$ is
equipped with the natural order and 
is a quotient of $[0,1]\setminus F$).

A verification shows that 
the above defines a map of sets
$$
|{X_\bullet}|=\varinjlim\limits_{{ F\subset [0,1] \text{ finite}}} {X_\bullet}(\pi_0([0,1]\setminus F)) 
\lra \homm{sSets}{\homm{\text{preorders}}{-}{[0,1]_\leq}}{ {X_\bullet} } 
$$

Conversely, a map $\pi:\homm{\text{preorders}}- {[0,1]_\leq}\lra {X_\bullet}$ of ssets determines a system of points  as follows:
$$\pi\left(\xra{\theta\,:\,n_\leq\lra[0,1]}\right)\in {X_\bullet}\left(\pi_0([0,1]\setminus\{\theta(0),..,\theta(n-1)\})\right)$$
and thereby a point of $|{X_\bullet}|$.

\subsubsection{The metric on the geometric realisation}
We make the following remark motivated by the definition of 
the metric $d_\mu(u,v)$, $u,v\in |X_\bullet|$, in [Drinfeld, following (1.7)]. 

Recall [Drinfeld] describe the topology on $|X_\bullet|$ using the following metric. 
We quote:
\bi\item
To define it first choose a measure $\mu$
on $I$ such that the measure of every point is zero and the
measure of every non-empty open set is non-zero. For a
finite $F\subset I$ denote by $\mu_F$ the measure on 
$\pi_0(I\setminus F)$ induced by $\mu$. If
$u,v\in X (\pi_0(I\setminus F))$ define the distance 
$d_{\mu}(u,v)$ to be the minimum of 
$\mu_F (\pi_0(I\setminus F)\setminus A)$ for all subsets
$A\subset\pi_0(I\setminus F)$ such that the images of
$u$ and $v$ in $X(A)$ are equal. If $F'\supset F$ and
$u',v'$ are the images of $u,v$ in $X (\pi_0(I\setminus F))$
then $d_{\mu}(u',v')=d_{\mu}(u,v)$, so we get a well defined
metric $d_{\mu}$ on $|X|$. 
\ei

For us it is more convenient to modify the definition slightly and
for 
$u,v\in X (\pi_0(I\setminus F))$ define the distance 
$d'(u,v)$ to be the minimum of 
$$\dist(F,A):=\max_{x\in F}\min_{a\in A}\dist(x,a)$$
for all subsets
$A\subset\pi_0(I\setminus F)$ such that the images of
$u$ and $v$ in $X(A)$ are equal.

Our notion of   $\theta\delta$-Skorokhod neighbourhood allows us to define a metric as follows:
take the largest metric satisfying the inequalities
\bi\item
$d''(x[1],x[2])\leq \delta$ whenever $1$-dimensional simplex $x$ lies in the  $\theta\delta$-Skorokhod neighbourhood for each $\theta$. 
\ei
Note that the uniform structure given by the metric can be defined entirely in terms 
of the neighbourhood structure on the Skorokhod space. Namely, this uniform structure 
is generated by entourages, for $\varepsilon$ a neighbourhood of $1$-simplicies 
\bi
\item 
$\left\{(\gamma_1[0],\gamma_2[0]) : \gamma_1, \gamma_2 \in \varepsilon,  \gamma_1[1]=\gamma_2[1]\right\} \subset |X_\bullet|\times |X_\bullet|$ 
\ei

The following construction relates the metrics $d'$ and $d''$. 

For $t_0\in [0,1]$ and a subset $A$, define a $A$-step function $\text{step}_{A,t_0}$: 
$t\longmapsto t$ if $t\leq t_0$, and $t\longmapsto \min\{ a\in A\cup\{1\} \,:\, a\geq t \} $.

For $u\in X_\bullet(\pi_0([0,1]\setminus F))$, consider the morphism 
$u_A: \hom-{2_\leq} \times [0,1]_\leq \lra X_\bullet$ defined as follows:
$$(1\leq l_1\leq ... \leq l_n\leq 2), (0\leq t_1\leq ...\leq t_n\leq 1) 
\longmapsto\,\,\,\,\,\,\,\,\,\,\,\,\,\,\,\,\,\,\,\,\,\,\,\,\,\,\,\,\,\,\,\,\,\,\,\,\,\,\,\,\,\,\,\,\,\,\,\,\,\,\,\,\,\,\,\,\,\,\,\,\,\,\,\,\,\,\,\,\,\,$$
$$\text{ }\longmapsto u\left[n_\leq\longmapsto [0,1]_\leq  
\xra{\text{step}_{A,\max\{t_i:\,l_i=1\}}} [0,1]_\leq \lra \pi_0([0,1]\setminus F) \right]$$
This morphism belongs to the $\theta\delta$-Skorokhod neighbourhood for any $\theta:n_\leq \lra 2_\leq$
and any $\delta>\dist(F,A)$ where $\dist(F,A):=\max_{x\in F}\min_{a\in A}\dist(x,a)$.
To see this, in the notation of the definition of the Skorokhod space,
for an $x= (t_1\leq ...\leq t_n)\in ([0,1]_\leq)_n$ pick  
$x'=(t_1\leq ...\leq t_n\leq a\leq a\leq ... \leq a)\in ([0,1]_\leq)_N$ where $a\in A\cup\{1\}$;
then $f(\theta,x'[N-n+1..N])=u_A(\theta,  a\leq a\leq ... \leq a)$
which is the value of the morphism $[0,1]_\leq\lra X_\bullet$ corresponding to $u\in X_\bullet(\pi_0([0,1]\setminus F))$,
and hence lies on the main diagonal in $(X_\bullet)_n$.

Now we may define the metric $d''':|X_\bullet|\times |X_\bullet|\lra \Bbb R $ to be the largest metric such that
$d(u,u_A)\leq \dist(F,A)$ for any 
 $u\in X_\bullet(\pi_0([0,1]\setminus F))$ as above.

A verification shows that all these metrics determines the same topology as the metric defined by [Drinfeld].

\subsubsection{Quillen edgewise subdivision} Following  [Grayson,\S3.1, esp.~Def.3.1.4, Def.3.1.8], 
define the following endofunctor on $\Dop$:
$$n\longmapsto 2n, \ f:m\ra n \longmapsto \{ n+i\mapsto n+f(i),\ n-i \mapsto n-f(i) , \text{ for }i=0,...,n-1\}$$
We can then define an endofunctor of $\sFilth$ by $X_\bullet\longmapsto X_\bullet\circ e$,
and consider various subdivisions of the geometric realisation of a simplicial set 
$$\!\!\!\!\!\!\!
\Homm\sFilth { [0,1]_\leq} {X_\bullet\circ e}
\,\,\ \Homm\sFilth { [0,1]_\leq\circ e} {X_\bullet}
\,\,\ \Homm\sFilth { [0,1]_\leq\circ e} {X_\bullet\circ e}
\,\,\ \Homm\sFilth { [0,1]_\leq} {X_\bullet}\circ e 
$$

\section{Appendix. Informal considerations leading to our definitions}

Here is a verbose manner we spell out informal considerations leading to our definitions.

\subsection{The approach of Drinfeld} 

To construct geometric realisation, [Drinfeld,\S1] considers the set of connected components of the unit interval without finitely many points
and uses that this set carries a canonical linear order.
Our interpretation is based on the observation that it is equivalent to consider instead the monotone non-decreasing function
on the unit interval contracting the connected components. This would be a function from the unit interval to a finite linearly ordered set.

Let us now sketch the relationship between our construction and that of Drinfeld 
on the level of points, not discussing topology. 

Recall that [Drinfeld, (1.1)] defines the geometric realisation $|X_\bullet|$ of a simplicial set $X_\bullet$ as
$$                                                                                                                                                             
|{X_\bullet}|:=\varinjlim\limits_{{ F\subset [0,1] \text{ finite}}} {X_\bullet}(\pi_0([0,1]\setminus F))$$
where $F$ runs through the set of all finite subsets of $[0,1]$. The notation 
$X_\bullet(\pi_0([0,1]\setminus F))$ might require some explanation, as follows.
The set $\pi_0([0,1] \setminus F )$ of connected components of 
$[0,1] \setminus F $ inherits the linear order from $[0,1]$ as a quotient of $[0,1]\setminus F$, 
hence is canonically isomorphic to an object of $\Delta$; hence we are justified
in denoting by 
${X_\bullet}(\pi_0([0,1]\setminus F))$ the functor $X_\bullet$ applied to that object of $\Delta$. 

Denote by $ \homm{\text{preorders}}{-}{[0,1]_\leq}$ the simplicial set corepresented by the linear order $[0,1]_\leq$ 
in the category of preorders; one may want to think of it as a ``thick'' interval or infinite dimensional simplex $\Delta_{[0,1]}$.

We observe that a finite subset $F\subset [0,1]$ and an $x\in {X_\bullet}(\pi_0([0,1]\setminus F))$ determines a morphism of sSets 
$ \homm{\text{preorders}}{-}{[0,1]_\leq} \lra  {X_\bullet} $ as follows:
$$ \homm{\text{preorders}}{n_\leq}{[0,1]_\leq} \lra  {X_\bullet}(n_\leq) $$
$$\overrightarrow t \longmapsto x\left[\xra{n_\leq\xra{\overrightarrow t} [0,1]\lra \pi_0([0,1]\setminus F)}\right]\in {X_\bullet}(n_\leq)
$$
where $[0,1]\lra \pi_0([0,1]\setminus F)$ is the obvious map contracting the connected components (we need to make
a convention where to send points of $F$).

A verification shows that 
the above defines a map of sets
$$
|{X_\bullet}|=\varinjlim\limits_{{ F\subset [0,1] \text{ finite}}} {X_\bullet}(\pi_0([0,1]\setminus F)) 
\lra \homm{sSets}{\homm{\text{preorders}}{-}{[0,1]_\leq}}{ {X_\bullet} } 
$$
from the points of the geometric realisation to the set $ \homm{sSets}{\homm{\text{preorders}}{-}{[0,1]_\leq}}{ {X_\bullet} }$
of morphisms in $\sSets$ from  ${\homm{\text{preorders}}{-}{[0,1]_\leq}}$ to  ${X_\bullet}$.

Conversely, a map $\pi:\homm{\text{preorders}}- {[0,1]_\leq}\lra {X_\bullet}$ of ssets determines a system of points  as follows:
$$\pi\left(\xra{\theta\,:\,n_\leq\lra[0,1]}\right)\in {X_\bullet}\left(\pi_0([0,1]\setminus\{\theta(0),..,\theta(n-1)\})\right)$$
and thereby a point of $|{X_\bullet}|$.

Note that our reformulation makes it immediate that, on the level of points,
the geometric realisation commutes with (arbitrary) injective limits, e.g. with Cartesian
products, and that the geometric realisation is equipped with an action of
orientation preserving homeomorphisms of [$0,1]$.

\subsection{Remark 2.4.1-2 of Grayson on the standard geometric simplex as a space of functions}

Let us now explain how to see construction as an attempt to reformulate in terms of $\sFilth$ 
a remark [Grayson, 2.4.1-2], cf.\S\ref{grayson:remark}, that the standard geometric simplex $\Delta_N\subset \RR^N$ 
can be viewed
as the space  of semi-continuous monotone functions $[0,1]_\leq \lra (N+1)_\leq$
with a Skorokhod-type metric defined as follows.

Following [Kolmogorov, \S2, Def.1], 
say that 
 two functions $f,g:[0,1]\lra \{1,..,N,N+1\}$ are {\em $\epsilon$-equivalent} or {\em at Skorokhod distance at most $\epsilon$} 
iff  
%
\bi\item[$(\star)$]
for each $t\in[0,1]$ there is $s\in [0,1]$ such that 
$$|t-s|<\epsilon\text{ and }\dist(f(t),g(s))<\epsilon$$ 
\ei 
and dually, for each $s\in[0,1]$ there is $t\in [0,1]$ 
$|t-s|<\epsilon\text{ and }\dist(f(t),g(s))<\epsilon$. 
Here we consider $\{1,..,N,N+1\}$ equipped with the discrete metric $\dist(m,n)=|m-n|$,
hence the condition $\dist(f(t),g(s))<\epsilon$ above is equivalent to $f(t)=g(s)$ for $\epsilon<1$.

This metric naturally occurs in probability theory  [Kolmogorov, \S2, Def.1]; 
the original goal of the definition was to define a distance or convergence
for (distributions of) stochastic processes such that
a small distortion of either timings of events
or their values results in a small distance.

\subsubsection{Rewriting Remark 2.4.1-2 of Grayson in simplicial language}
First use the Yoneda lemma to rewrite $\homm{\text{preorders}}{[0,1]_\leq}{N_\leq }$ in terms 
of simplicial sets as 
$$\homm{\sSets}{\homm{\text{preorders}}-{[0,1]_\leq}}{\homm{\text{preorders}}-{N_\leq}}$$

Now we may say that in $(\star)$, assuming that $t\leq s$ and $f(t)\leq g(s)$, 
that the pair $(f(t)\leq g(s))$ is a $1$-simplex of $\Delta_N$,
that the pair $(t\leq s)$ is a $1$-simplex of $\Delta_{[0,1]_\leq}:=\homm{\text{preorders}}-{[0,1]_\leq}$,
and that $0$-simplex $t$ is a face of $(t\leq s)$. A pair of function $f,g:[0,1]\lra N_\leq$
is a $1$-simplex of the inner hom 
$$\Homm{\sSets}{\homm{\text{preorders}}-{[0,1]_\leq}}{\homm{\text{preorders}}-{N_\leq}},$$
again provided that $f(t)\leq g(t)$ for $t\in [0,1]$, or, equivalently, 
 $f(t)\leq g(s)$ whenever $t\leq s$. 

Hence, we rewrite $(\star)$ as 
\bi\item[$(\star')$]
for each $0$-simplex $t$ of $\Delta_{[0,1]}$ 
there is $1$-simplex $(t\leq s)$ of $\Delta_{[0,1]_\leq}$ such that 
\bi \item $t$ is a face of $(t\leq s)$
\item both $1$-simplex $(t\leq s)$ and its image $(f(t),g(s))$ are ''$\epsilon$-small''.
\ei 
\ei
Now we need to give an exact meaning to the phrase ''$\epsilon$-small''
without any reference to real numbers.

\subsubsection{
Introducing the additional structure of topological flavour on simplicial sets: generalities}
The notion of smallness comes from the topological and/or metric structure on $[0,1]$,
hence we should endow the sset $\Delta_{[0,1]_\leq}=\homm{\text{preorders}}-{[0,1]_\leq}$ with 
an additional structure of topological flavour; we do it as follows.


Recall that [Bourbaki,Introduction] ``a topological structure on a set enables one
to 
give an exact meaning to the phrase
``whenever $x$ is sufficiently near $a$, $x$ has the property $P\{x\}$''
by saying that the elements with property $P\{x\}$ form a neighbourhood of $x$.  
Hence, we may say that ``a set $E$ carries a topological structure whenever we have 
associated with each element of $E$, by some means or other, a family 
of subsets of $E$ which are called neighbourhoods of this element --- provided 
of course that these neighbourhoods satisfy certain conditions (the axioms 
of topological structures).''

In a similar way, we equip a simplicial set 
with an additional structure which enables one
to give an exact meaning to the phrase
``whenever a simplex $x$ is sufficiently small,  $x$ has the property $P\{x\}$''; 
we say that a simplicial set $E_\bullet$ carries this additional structure  
whenever we have 
associated with each dimension $n\geq 0$, by some means or other, a family 
of subsets of $E_n$ which are called {\em neighbourhoods} --- provided 
of course that these neighbourhoods satisfy certain conditions.

These conditions in fact say that 
the simplicial set with this additional structure is  a simplicial object 
in the category of filters in the sense of Bourbaki
defined in \S\ref{filter:def}.

\subsubsection{The additional structure in our case}
Let us now describe explicitly the sset $\Delta_{[0,1]_\leq}$ and its extra structure ``inherited'' from $[0,1]$ : 
an $n$-simplex is a non-decreasing sequence $(t_0\leq ... \leq t_n)$ where $ t_0,t_1,..,t_n\in [0,1]$,
and a subset of the set of $n$-simplicies is {\em big} or {\em a neighbourhood} iff 
it contains for some $\epsilon>0$ the ``$\epsilon$-neighbourhood of the main diagonal''
$$
\{\,(t_0\leq ... \leq t_n)\,:\, t_n\leq t_0+\epsilon,\, t_0,..,t_n\in [0,1]\}
$$

Thus, in $(\star')$, assuming as before $t\leq s$, we may express $|t-s|<\epsilon$
by saying that the $1$-simplex $(t\leq s)$ belongs to the $\epsilon$-neighbourhood of the main diagonal.

In $\Delta_N$, an $n$-simplex is a non-decreasing tuple $(t_0\leq ... \leq t_n)$ where $t_0,..,t_n\in \{1,..,N+1\}$;
a subset of the set of $n$-simplicies is {\em big} or {\em a neighbourhood} iff 
it contains the ``the main diagonal''
$$
\{\,(t,t,..t)\,:\, t\in\{0,1,..,N+1\}\,\}
$$

Thus, in $(\star')$, assuming $f(t)\leq g(s)$ and $\epsilon<1$, we may express $\dist(f(t),g(s))<\epsilon$, i.e.~$f(t)=g(s)$,
by saying that the $1$-simplex $(f(t)\leq g(s))$ belongs to neighbourhood formed by the main diagonal (and therefore 
to each neighbourhood of this structure).

\subsubsection{Neighbourhoods in Skorokhod metric in simplicial language} 

Now we are ready to turn $(\star)$ into a definition of  {\em Skorokhod neighbourhoods} in the set of $1$-simplicies of 
the inner hom $$\Homm{\sSets}{\homm{\text{preorders}}-{[0,1]_\leq}}{\homm{\text{preorders}}-{N_\leq}}\,\,:$$
given a neighbourhood $\varepsilon=\{(t\leq s): t\leq s\leq t+\epsilon,\, t,s\in [0,1]\}\subset\Delta_{[0,1]}(1)$ 
and a neighbourhood $\delta=\{(n,n):n=1,..,N+1\}\subset \Delta_{N}(1)$, {\em the Skorokhod $\varepsilon\delta$-neighbourhood} is
the set of all $1$-simplicies $(f\leq g)$ 
\bi\item[$(\star'')$]
for each $0$-simplex $t$ of $\Delta_{[0,1]}$ 
there is $1$-simplex $(t\leq s)$ of $\Delta_{[0,1]_\leq}$ such that 
\bi \item $t$ is a face of $(t\leq s)$
\item $1$-simplex  $(t\leq s)\in\delta$ and its image $(f(t),g(s))\in\varepsilon$
\ei 
\ei

If we continue similar analysis, we arrive at the definition in \S\ref{skoro:filter} of the Skorokhod filter
and in  \S\ref{skoro:space}
of the Skorokhod space $\HommS\Filt{X_\bullet}{Y_\bullet}$ for arbitrary objects $X_\bullet,Y_\bullet$ of $\sFilth$.

\section{\label{remarks}Appendix. An assortie of examples}

We include some examples to indicate the expressive power of $\sFilth$, or rather
to indicate a way how to use $\sFilth$ to rephrase certain standard familiar definitions. 
Their verification is routine
and a matter of working through the notation and definitions.   

This section reports on an unfinished work unlikely to be continued; we do apologise. 
Some details can be found in
\href{http://mishap.sdf.org/6a6ywke/6a6ywke.pdf}{[6]}. Notably, \href{http://mishap.sdf.org/6a6ywke/6a6ywke.pdf}{[6,\S3]} 
includes a discussion of the 
intuition of point-set topology in $\sFilth$.

\subsection{Compactness and compact-open topology} A topological space is compact iff every open covering has a finite subcovering. 
We now introduce enough notation to rewrite this definition in a way which makes obvious how to generalise it to $\sFilth$.

\subsubsection{Transcribing ``every covering has a finite subcovering''}
View a {\em covering $U_x\ni x, x\in X$} of a topological space $X$ as a subset of $X\times X$
$$\varepsilon:=\bigcup\limits_{x\in X} \{x\}\times U_x \subset X\times X$$
Viewed this way, coverings of $X$ form a filter {\em of coverings} on $X\times X$; we call its elements {\em neighbourhoods} 
and denote them by $\varepsilon,\delta,...$. 
Instead of saying {\em point $y$ lies in the neighbourhood $U_x$ of point $x$} we say
that {\em the simplex $(x,y)$ is $\varepsilon$-small} meaning that $(x,y)\in \varepsilon\subset X\times X$.

The covering $U_x$ has a finite subcovering iff $\emptyset$ is contained in the filter generated by the complements $\bar U_x:=X\setminus U_x$ 
of $U_x,x\in X$. Denote by $$\bar\varepsilon_x := 
\{ y\,:\, (x,y) \text{ is not }\varepsilon\text{-small}\,\}=\{ y\,:\, (x,y) \not\in \varepsilon\}=X\setminus U_x$$ 
the ``tail'' fibres over $X$ of the complement $\bar\varepsilon:=X\times X\setminus\varepsilon$ of neighbourhood $\varepsilon$.
Denote by $ \left< \bar\varepsilon_x\right>_{ x\in X} $ the filter on $X$ generated by the fibres of the complement $\bar\varepsilon$.
Then we may that
\bi\item a topological space $X$ is compact iff for each neighbourhood $\varepsilon\subset X\times X$ in the filter of coverings
$$
\emptyset \in \left< \bar\varepsilon_x\right>_{ x\in X}
$$
\ei
\subsubsection{Quantifying over a neighbourhood} 
It is somewhat unnatural, either in topology but especially in the context of $\sFilth$,
 to write  ``for all $x$ in $X$'' rather than ''for all $x$ in a neighbourhood $\delta$'', hence we rewrite this further 
and consider instead $\emptyset \in \left< \bar\varepsilon_x\right>_{ x\in \delta}$. However, this is unnecessary in this case,
as we may replace $\varepsilon$ by the intersection with the preimage of $\delta$:
$$
\left< \bar\varepsilon_x\right>_{ x\in \delta} 
=
\left< \overline{(\varepsilon\cap \delta^{-1})}_x\right>_{ x\in X}
$$

\subsubsection{Generalising to arbitrary dimension and objects}
To generalise this to arbitrary dimension and an arbitrary object of $\sFilth$, we only need to make the following convention:
for $n>m>0$, a point $x\in X_\bullet(m_\leq)$, and a neighbourhood $\varepsilon\subset X_\bullet(n_\leq)$,
denote by $\varepsilon_{x:[1..m]}$  the fibre over $x$ of the face map $[1..m]:X_\bullet(n_\leq)\lra X_\bullet(m_\leq)$
to the first few coordinates.

\subsubsection{A definition of compactness and being proper}

Call an object $X_\bullet:\sFilth$ {\em compact} iff  
\bi\item 
iff for each $n>m>0$, for each neighbourhood $\varepsilon\subset X_\bullet(n_\leq)$ 
$$
\emptyset \in \left< \bar\varepsilon_{x:[1..m]}\right>_{ x\in X}
$$
\ei

Call a morphism $f:X_\bullet\lra Y_\bullet$ in $\sFilth$ {\em proper} iff
\bi\item for each $n>m>0$,
for each neighbourhood $\varepsilon$ in $X_\bullet(n_\leq)$ 
there exists a neighbourhood $\delta$ in $Y_\bullet(m_\leq)$ 
for each $y_0\in X_\bullet(m_\leq)$
$$
\emptyset\in \left< \bar\varepsilon_{x:[1..m]} \cap f^{-1}(\delta) \right>_{ f(x)=y_0[1..m]} 
$$
%
\ei

Let us rewrite the latter back in the notation of coverings, for $X_\bullet$ and $Y_\bullet$ 
objects associated with topological spaces, $n=2>m=1>0$:
\bi
\item given a covering of $X$, each point of $y$ has a small enough neighbourhood such that its preimage has a finite subcovering
\item (the same in notation)
for each covering $(U_x)_{x\in X}$ (i.e.~a neighbourhood $\varepsilon:=\cup_{x\in X} \{x\}\times U_x$ in the filter of coverings)
for each point $y\in Y$ there is a neighbourhood $V_y\ni y$ of $y$ (i.e.~a neighbourhood $\delta:=\cup_{y\in Y} \{y\}\times V_y$ in the filter of coverings)
such that the preimage of $V_y$ is covered by a finite subcovering of $(U_x)_{x\in X}$. 
\ei
Thus it is the usual definition of a proper map stated in terms of open coverings.


\subsubsection{Core compact and compact-open topology}
Recall a subset $K\subset X$ of a topological space is called {\em core compact (in $X$)} iff every covering of $X$ has a finite subcovering of $K$. 
In $\sFilth$, we say that $K\subset X_\bullet(n_\leq)$ is {\em core compact (in $X_\bullet$)} iff 
$$
X_\bullet(n_\leq)\setminus K \in \left< \bar\varepsilon_{x:[1..m]}\right>_{ x\in X}
$$
 for each $n>m>0$, for each neighbourhood $\varepsilon\subset X_\bullet(n_\leq)$.

Note that the image of a core compact subset is necessarily a core compact subset, and, similarly to topology,
This leads to a functorial {\em core compact-open} filter structure on $\homm{\sFilth} {X_\bullet} {Y_\bullet}$
generated by neighbourhoods 
$$
\delta_{K,\varepsilon}:=\left\{ f: X_\bullet\lra Y_\bullet \,\,:\,\,  f(K)\subset \varepsilon\,\,\right\}
$$
where   $K\subset X_\bullet(n_\leq)$ is core compact, and $\varepsilon\subset X_\bullet(n_\leq)$ is a neighbourhood,
for some $n>0$. 

This functorial filter structure on Homs in $\sFilth$ would allow to define a functor 
$\sFilth^{\text{op}}\times \sFilth \lra \sFilth$
analogous to the notion of the space of maps with core-compact-open topology.
Let us spell some details.

Equip  $N$-simplicies $\hom{-}{N_\leq}$ with the antidiscrete filters, and equip 
each $\hom{\hom{-}{N_\leq}\times {X_\bullet}}{{Y_\bullet}}$ with the resulting core-compact-small filters.
Now consider the inner hom in $\sSets$ equipped with these filters, and 
its subsset  consisting of simplicies all of whose $0$-dimensional faces 
are continuous functions, i.e. belong to $\homm\sFilth{X_\bullet}{Y_\bullet}$. 
This defines a functor 
$$
\Homm\sFilth - - : \sFilth^\text{op}\times \sFilth \lra \sFilth 
$$ 
which ``enriches'' the $\text{Hom}$-functor in $\sFilth$.

\subsection{The bounded-small space of maps} 
Another analogue of the definition of the compact-open topology on the space of maps
may be formulated as follows.

Informally, call a simplex in a simplicial filter{\em Archimedean} iff it can
``split'' into a simplex with
arbitrarily small consecutive faces; call a set of Archimedean simplicies {\em bounded} 
there is a bound on the dimension (number of faces) of the simplicies with arbitrary small consecutive faces
one needs to consider. Intuitively, we'd like to think of such a simplex with arbitrarily small 
consecutive faces as an approximate homotopy. 

This is made formal as follows.

Call a subset $B\subset X_k$ {\em bounded} iff 
for each $n,m\geq 0$ for each neighbourhood $\varepsilon\subset X_n$
there is $N\geq 0$ such that for every $b\in B$ there exists $b_N\in X_N$ such that
\bi
\item $b$ is a face of $b_N$
\item the face $b_N[l_1\leq ... \leq l_n]\in \varepsilon$ is $\varepsilon$-small whenever $1\leq l_1\leq ...\leq  l_n\leq l_1+m\leq N$.
\ei

An alternative (possibly non-equivalent) definition is: call a subset $B\subset X_k$ {\em bounded} iff 
for each $A_\bullet,Y_\bullet:\sFilth$, each morphism $f: A_\bullet\times X_\bullet\lra Y_\bullet$, each $n\geq 0$, and each neighbourhood $\varepsilon\subset Y_n$ 
there exists a neighbourhood $\alpha\subset X_n$ such that $f(\alpha\times B)\subset \varepsilon$. 

With either definition, bounded sets form an ideal, i.e.~are closed under finite unions and taking subsets;
the image of a bounded set under a continuous map is necessarily bounded. 
Define the {\em bounded-small} filter on $\homm\sFilth {X_\bullet} {Y_\bullet}$ as follows:
a subset $\delta\subset \homm\sFilth {X_\bullet} {Y_\bullet}$ is a neighbourhood iff 
there is $n\geq 0$, a bounded set $B\subset X_n$ and a neighbourhood $\varepsilon\subset Y_n$
such that $f\in\delta$ whenever 
\bi
\item $f(b)\in \varepsilon$ whenever $b\in B$
\ei

Equip  $N$-simplicies $\hom{-}{N_\leq}$ with the antidiscrete filters, and equip 
each $\hom{\hom{-}{N_\leq}\times {X_\bullet}}{{Y_\bullet}}$ with the resulting bounded-small filters.
Now consider the inner hom in $\sSets$ equipped with these filters, and 
its subsset  consisting of simplicies all of whose $0$-dimensional faces 
are continuous functions. 

This defines a functor 
$$
\Homm\sFilth - - : \sFilth^\text{op}\times \sFilth \lra \sFilth 
$$ 
which ``enriches'' the $\text{Hom}$-functor in $\sFilth$. 

Perhaps more interesting is to consider the subfunctor of the inner hom in $\sSets$
formed by simplicies Archimedean with respect to the bounded-small filters.
$$
\hHomm\sFilth - - : \sFilth^\text{op}\times \sFilth \lra \sFilth 
$$

\subsection{Topological spaces and uniform structures as simplicial filters} 
  
A uniform structure {[Bourbaki, II\S1.1,Def.1]} on a set $X$ 
gives rise to a simplicial filter on the simplicial set  co-represented by the set $X$:
the filter on $X$ is antidiscrete, the filter on $X\times X$ is the uniform structure (i.e.~the filter of entourages), 
and the filter on $X^n,n>2$ is the coarsest filter such that all the face maps $X^n\lra X^2$ are continuous. 
Explicitly, in terms of a metric on $X$, the filter on $X^n,n>0$ is the filter of ``$\epsilon$-neighbourhoods 
of the 
main diagonal $\{ (x,x,...,x) : x\in X\}$'', i.e. the filter generated by sets, for $\epsilon>0$, 
$$\{(x_1,...,x_n)\in X^n: \dist(x_i,x_j)<\epsilon\}.$$ 
A $1$-simplex $b=(x,y)$, $x,y\in M$, is Archimedean iff for each $\varepsilon>0$ 
``there is an $\varepsilon$-approximate homotopy from $x$ to $y$'', i.e.~the points 
$x$ and $y$ can be joined by a sequence of points $b'_N = (x=x_1,x_2...,x_{N}=y$
at distance less than $\varepsilon$, $x_i,x_{i+1}<\varepsilon$ for all $1\leq i<i+1\leq N$. 
Two maps $f,g:M_1\lra M_2$ of sufficiently nice metric spaces (e.g.~for each 
$\varepsilon$-ball there is $\varepsilon'>0$ such any $\varepsilon$-ball contains
a contractible subset containing $\varepsilon'$-ball around its centre)
 are homotopic iff they determine an Archimedean $1$-simplex
in the bounded-small space of maps: in the notation of previous subsection, 
an $\varepsilon$-approximate homotopy is given by faces
$f=b_N[1],b_N[2],...,b_N[N]=g:M_1\lra M_2$ of $b=(f,g):M_1\times M_1\lra M_2\times M_2$.

In a similar way a topology on a set $X$ 
gives rise to a simplicial filter on the simplicial set  co-represented by the set $X$:
as before, the filter on $X$ is antidiscrete, and 
the filter on $X^n,n>2$ is the coarsest filter such that all the face maps $X^n\lra X^2$ are continuous,
and on $X\times X$ take the filter of ``non-uniform'' neighbourhoods of the main diagonal, i.e.~sets of form
$$\bigcup\limits_{x\in X, \text{ and }U_x\text{ is a neighbourhood of }x} \{x\}\times U_x $$

In this way, the category of topological spaces and that of uniform spaces 
fully faithfully embed into $\sFilth$.

\subsection{\label{quasi:iso}Quasi-isomorphisms in large scale geometry} Given a metric on a set $X$, consider the filter on $X^n$ generated by
$$\{ (x_1,...,x_n): \text{ for all }i,j\text{ either }\dist(x_i,x_j)>D \text{ or } x_i=x_j\,\}$$ this associates  an object of $\sFilth$ 
with a metric space $X$. This filters ``remembers'' some of the structure considered in large scale geometry [Gromov,\S0.2.BC],
for example two quasi-geodesic spaces are quasi-isomorphic iff their associated 
objects are isomorphic in $\sFilth$.

\subsection{\label{fibre:bundle}Locally trivial fibre bundles} Let $[+1]:\Dop\lra\Dop$ be the endomorphism of $\Dop$ adding a new minimal element, 
i.e.~$n[+1]:=n+1$, $n\in \NN$,  and for a morphism $f$,  $f[+1](0):=0$,  $f[+1](i+1):=f(i)+1$, $0\leq i<n$. 
A morphism $p:X\lra B$ of topological spaces is locally trivial with fibre $F$ iff in $\sFilth$ 
the map $p_\bullet:X_\bullet\lra B_\bullet$ of their associated objects becomes ``trivial'' after 
pullback along the obvious map $B_\bullet\circ[+1]\lra B_\bullet$ ``forgetting the first coordinate'': 
there is an isomorphism in $\sFilth$ over $B_\bullet\circ [+1]$
$$
(B_\bullet\circ [+1]) \times_{B_\bullet} X_\bullet \xleftrightarrow{(iso)} (B_\bullet\circ[+1]) \times F_\bullet 
\,\,\,\,\,\,\,\text{ over } B_\bullet\circ [+1]$$

\subsection{\label{path:space}A notion of path and loop space}
Let $I_\leq$ be a linear order, 
and let $\Bbb I_\leq$ denote 
the simplicial set 
 co-represented by $[0,1]_\leq$ in the category of preorders 
$$ n_\leq\longmapsto \homm{{\text{preorders}}}{n_\leq}{[0,1]_\leq}$$ 
equipped with the the subdivision filters, cf.~\ref{unit-interval}. 
 
Let  $X_\bullet:\sFilth$ be a simplicial filter. The {\em $I_\leq$-path space of $X_\bullet$}
is 
$$\hHomm{\sFilth}{\Bbb I_\leq} {X_\bullet} $$

Now let $[+1]:I_\leq\lra I_\leq$ be an endomorphism of $I_\leq$ defining an action $\NN\times I_\leq \lra I_\leq$
of $\NN$ on $I_\leq$.

 {\em A  $I_\leq/\NN$-loop space} of a simplicial filter $X_\bullet:\sFilth$ is 
the subsspace of the  $I_\leq$-path space of $X_\bullet$ consisting of simplicies
Archimedean in the subsspace of simplicies invariant under the $\NN$-action.

More explicitly this can be described as follows.

\subsection{\label{path:space:explicit} Path and loop spaces: an explicit description}
Consider the simplicial set 
$$\Homm{\sFilth}{\Homm{\text{preorders}}-{ I_\leq}} {X_\bullet} .$$
Let $m>0$ and $\varepsilon\subset X(n_\leq)$ be a neighbourhood.

We say that {\em $\gamma\in \hom{\homm{\text{preorders}}-{I_\leq}} {X_\bullet}(n_\leq) $
is an $\varepsilon$-fine discretised homotopy at a mesh $\overline t:N\lra n$} whenever for each $\theta:N\lra n$
$\gamma(\theta,\overline t)[i+1..i+m]$ is $\varepsilon$-small, i.e.~$\gamma(\theta,\overline t)[i+1..i+m]\in\varepsilon$.

We say that $\overline {t'}:N_\leq \lra I_\leq $ is  {\em a refinement of $\overline t:m_\leq\lra I_\leq$} iff 
for some $1\leq i_1\leq .. \leq i_m\leq N$ it holds $\overline t=\overline{t'}[i_1\leq .. \leq i_m]$, 
i.e.~the simplex $\overline t$ is a face of simplex $\overline{t'}$. 

We say that  {\em $\gamma\in \hom{\homm{\text{preorders}}-{I_\leq}} {X_\bullet}(n_\leq) $ 
is a homotopy} iff for each $k,m>0$, each $\overline t:k_\leq\lra I_\leq$, 
and each neighbourhood $\varepsilon\subset X_m$  there is $N>0$ and a refinement   
$\overline{t'}:N_\leq \lra I_\leq $ such that  $\gamma\in \hom{\homm{\text{preorders}}-{I_\leq}} {X_\bullet}(n_\leq) $ 
is an $\varepsilon$-fine discretised homotopy at the mesh $\overline{t'}:N\lra n$.

The {\em $I_\leq$-path space} of a simplicial filter $X_\bullet$
is the subsset of $$\hom{\homm{\text{preorders}}-{I_\leq}} {X_\bullet} $$
formed by those simplicies which are homotopies and equipped with filters as follows. 
A subset of $\hom{\homm{\text{preorders}}-{I_\leq}} {X_\bullet}(n_\leq) $ 
is a neighbourhood iff it contains all the homotopies $\varepsilon$-finely discretised at any refinement of $\overline t:N\lra n$, 
for some neighbourhood $\varepsilon\subset X_m$ and some increasing sequence $t_1\leq ... \leq t_N$ of elements of $I_\leq$.

Now let $[+1]:I_\leq\lra I_\leq$ be an endomorphism of $I_\leq$ defining an action $\NN\times I_\leq \lra I_\leq$
of $\NN$ on $I_\leq$.
The definition of {\em a $I_\leq/\NN$-loop} in a simplicial filter $X_\bullet:\sFilth$ is analogous to that of a $I$-path
but we require ``meshes'' $t:k_\leq \lra I$ and $t':N_\leq \lra I $ to be invariant under the action of $\ZZ$.

We say that {\em $\gamma\in \hom{\homm{\text{preorders}}-{I_\leq}} {X_\bullet}(n_\leq) $
is an $\varepsilon$-fine discretised homotopy of loops at a mesh $\overline t:N\lra n$} whenever
\bi\item for each $N',k>0$, each $\theta:N'\lra n$, each $\overline t': N'_\leq\lra I_\leq $ it holds 
$\gamma(\theta,\overline {t'})=\gamma(\theta, \overline {t'}+k)$
\item for each $\theta:N\lra n$ and each $k>0$
$\gamma(\theta,\overline t+k)[i+1..i+m]$ is $\varepsilon$-small, i.e.~$\gamma(\theta,\overline t+k)[i+1..i+m]\in\varepsilon$.
\ei

We say that  {\em $\gamma\in \hom{\homm{\text{preorders}}-{I_\leq}} {X_\bullet}(n_\leq) $ 
is a homotopy of $I_\leq/\NN$-loops} iff for each $k,m>0$, each $\overline t:k_\leq\lra I_\leq$, 
and each neighbourhood $\varepsilon\subset X_m$  there is $N>0$ and a refinement   
$\overline{t'}:N_\leq \lra I_\leq $ such that  
$\gamma\in \hom{\homm{\text{preorders}}-{I_\leq}} {X_\bullet}(n_\leq) $ 
is an $\varepsilon$-fine discretised homotopy at the mesh $\overline{t'}:N\lra n$.

 {\em A  $I_\leq/\NN$-loop space} of a simplicial filter $X_\bullet:\sFilth$ is 
the subsspace of the  $I_\leq$-path space of $X_\bullet$ consisting of homotopies of $I_\leq/\NN$-loops.

\subsection{\label{ramsey}Ramsey theory and indiscernible sequences} Let $c$ be a colouring of {\em non-degenerate}
simplicies of dimension $n$ in a simplicial set $X_\bullet$. Call a simplex $c$-homogeneous iff
all of its non-degenerate faces of dimension $n$ have the same $c$-colour; trivially any $0$-simplex
is $c$-homogeneous for any $c$. Call a subset of $X_N$ a {\em $c$-neighbourhood of the main diagonal}
iff it contains all the $c$-homogeneous simplicies. 

Let $M$ be a model of a first-order theory, and let $M_\bullet=\left(n_\leq \mapsto M^n\right)$
 be the simplicial set co-represented by the set $M$ of its elements. Consider a formula $\phi(x_1,..,x_n)$ 
with $n$ variables as a colouring on the sequences of $n$ {\em distinct} elements of $M$; 
in terminology of model theory, 
a $\phi$-homogeneous non-degenerate simplex in $M_\bullet$ is called a $\phi$-indiscernible sequence 
of finite length (usually one considers indiscernible sequences whose length is a cardinal). 
Equip $M_\bullet$ with the filter of $\phi()$-neighbourhoods of the main diagonal
for all the formulas $\phi()$; this defines a space analogous to the Stone space of types 
of a model $M$ but capturing something about indiscernible sequences.

\subsection{Limits of sequences.}\label{limits-samples}
\def\NNcof{\NN_\text{cof}} \def\NNNcof{ \left( n\mapsto \NNcof^n\right)}
Let $\NNcof$ be the filter of cofinite subsets, and let $\NNNcof$ denote the simplicial object $\Dop\lra\Filth,\ n_\leq \mapsto \NNcof^{n}$ 
of Cartesian powers of $\NNcof$. Let $ M_\mU$ denote the object of $\sFilth$ corresponding to a metric (uniform) space as defined above. 

To give a {\em Cauchy sequence $(a_n)_{n\in\NN}$} in a metric space $M$ is
to give a morphism 
$\bar a:\NNNcof \lra M_\mU$, $(i_1,..,i_n)\mapsto (a_{i_1},...,a_{i_n})$.
Indeed, such a morphism is necessarily induced by a function $\NN\lra M$ of sets;
continuity in $\sFilth$ means that 
for each $\epsilon>0$ the preimage of $\varepsilon:=\{(x,y):\dist(x,y)<\epsilon\}$ contains  $\delta:=\{(n,m):n,m>N\}\}$
for some $N$ large enough, i.e.~$\dist(a_n,a_m)<\epsilon$ for $n,m>N$.
The sequence $(a_n)_{n\in\NN}$  {\em converges} 
iff the morphism $\bar a: \NNNcof \lra M_\mU$ 
factors as  $\bar a: \NNNcof \xra{\bar a_\infty} M_\mU[+1] \xra{[-1]}  M_\mU$.
Moreover, the
 morphism ${\bar a_\infty}:  \NNNcof \xra{\bar a_\infty} M_\mU[+1]$
is necessarily of form 
$M^n\lra M^{n+1}, (i_0,i_1,...,i_n)\mapsto (a_\infty, a_{i_1},...,a_{i_n})$
where $a_\infty$ is the {\em limit} of sequence $(a_n)_{n\in\NN}$.

To see this, first note that
the underlying sset of  $M_\mU[+1]$ is 
 a disjoint union $M_\mU[+1]=\sqcup_{a\in M} \{a\}\times M_\mU$
of copies of $M_\mU$, and that the underlying sset of $ \NNNcof$ is connected.
Hence, to pick a factorisation of the underlying ssets 
is to pick an $a\in M$. Now, the map $ \NNNcof \lra \{a\}\times M_\mU$
is continuous iff for each $\epsilon>0$ the preimage of $\varepsilon:=\{(a,x):\dist(a,x)<\epsilon\}$ 
contains  $\delta:=\{(n,m):n,m>N\}$, i.e.~for $m>N$ $\dist(a,a_m)<\epsilon$.

A {\em limit} of a sequence in a topological space is defined with help of factorisation
 $$\bar a:  \left( n\mapsto (\NNcof^n)_\diag\right) \xra{\bar a_\infty} X_\ttt[+1] \xra{[-1]}  X_\ttt$$
where $X_\ttt$ denotes the simplicial object corresponding to a topological space $X$, and 
 $ \left( n\mapsto (\NNcof^n)_\diag\right)$ denotes  the simplicial object $\Dop\lra\Filth,\ n_\leq \mapsto (\NNcof^{n})_\diag$
of Cartesian powers equipped with  ``the filter of cofinite diagonals'', i.e.~a subset of $\NN^n$ is a neighbourhood in $ (\NNcof^{n})_\diag$
iff it contains the set $\{ (i,i,..,i): i>N\}$ for some $N>0$.

\subsection{Compactness and completeness} The considerations above allow to express {\em completeness} as a lifting property 
in the sense of Quillen: a metric  space is complete, i.e.~each Cauchy sequence converges, iff $\bot\lra \NNNcof$ 
is orthogonal to $M_\mU[+1]\lra M_\mU$. This suggests one to consider a class of complete metric spaces defined as 
the double left-right orthogonal of the morphism $\RR_\mU[+1]\lra \RR_\mU$ associated with $\RR$, 
an archetypal example of a complete metric space. See \href{http://mishap.sdf.org/6a6ywke/6a6ywke.pdf}{[6,\S2.1.4,\S4.11]}
for speculations about compactness and completeness in terms of lifting properties and $\sFilth$.

\subsection{Uniform continuity and equicontinuity}

A {\em uniformly continuous function} $f:L\lra M$ is a morphism $L_\mU\lra M_\mU$. 
Indeed, for every  $\epsilon>0$ the preimage of $\varepsilon:=\{(u,v):\dist(u,v)<\epsilon\}\subset M\times M$ 
contains  ${{\updelta}} :=\{(x,y): \dist(x,y)<\delta\}\subset L\times L$ for some $\delta>0$, 
i.e.~$\dist(f(x),f(y))<\varepsilon$ whenever  $\dist(x,y)<\delta\}$.

A {\em uniformly equicontinuous sequence $(f_i:L\lra M)_{i\in\NN}$ of uniformly continuous functions} 
is a morphism $  \left( n\mapsto \NNcof\right) \times L_\mU \lra M_\mU$, 
$$\NN\times L^n\lra M^n,\ (i,x_1,...,x_n) \mapsto ( f_{i}(x_1),...,f_{i}(x_n)),$$
or, equivalently,  
is a morphism $  \left( n\mapsto (\NNcof^n)_\diag\right) \times L_\mU \lra M_\mU$, 
$$\NN^n\times L^n\lra M^n,\ (i_1,...,i_n,x_1,...,x_n) \mapsto ( f_{i_1}(x_1),...,f_{i_n}(x_n))$$
where  $ \left( n\mapsto (\NNcof^n)_\diag\right)$ denotes  the simplicial object $\Dop\lra\Filth,\ n_\leq \mapsto (\NNcof^{n})_\diag$
of Cartesian powers equipped with  ``the filter of cofinite diagonals'', i.e.~a subset of $\NN^n$ is a neighbourhood of $ (\NNcof^{n})_\diag$
iff it contains the set $\{ (i,i,..,i): i>N\}$ for some $N>0$.

Indeed, for every  $\epsilon>0$ the preimage of $\varepsilon:=\{(u,v):\dist(u,v)<\epsilon\}\subset M\times M$ 
contains  ${\updelta}:=\{(n,x,y):n>N,\dist(x,y)<\delta\}\subset \NN\times \NN\times L\times L $, 
resp.,  ${\updelta}:=\{(n,n,x,y):n>N,\dist(x,y)<\delta\}\subset \NN\times \NN\times L\times L $,
for some $\delta>0$ and $N>0$, 
i.e.~$\dist(f_n(x),f_n(y))<\varepsilon$ whenever $n>N$ and $\dist(x,y)<\delta$.

A {\em uniformly Cauchy sequence $(f_i:L\lra M)_{i\in\NN}$ of uniformly continuous functions} 
is a morphism $ \NNNcof \times L_\mU \lra M_\mU$, 
$$\NN^n\times L^n\lra M^n,\ (i_1,...,i_n,x_1,...,x_n) \mapsto ( f_{i_1}(x_1),...,f_{i_n}(x_n)).$$

Indeed, for every  $\epsilon>0$ the preimage of $\varepsilon:=\{(u,v):\dist(u,v)<\epsilon\}\subset M\times M$ 
contains  ${\updelta}:=\{(n,m,x,y):n,m>N,\dist(x,y)<\delta\}\subset \NN\times \NN\times L\times L $
for some $\delta>0$ and $N>0$, 
i.e.~$\dist(f_n(x),f_m(y))<\varepsilon$ whenever $n,m>N$ and $\dist(x,y)<\delta$.

Convergence and limits of sequences of functions can be reformulated in a
manner similar to limits of sequences of points above,
but somewhat more complicated, see \href{http://mishap.sdf.org/6a6ywke/6a6ywke.pdf}{[6,\S2.1.3,\S4.10]}.

\subsection{Open problems} It is premature to formulate any open problems and conjectures
before it has been verified whether the definitions above have the expected properties
and admit development of a theory. Nevertheless, we suggest a couple of suggestions 
for further research.

The reformulation of the notion of locally trivial fibration in \S\ref{fibre:bundle} 
suggests the following direction. Note that in the $\sFilth$ formalism does not suggest 
how to reverse direction of a homotopy $f\rightsquigarrow g$ to get $g \rightsquigarrow f$.

\subsubsection{Axiomatize homotopy theory} Formulate ``non-commutative'' axioms of homotopy theory 
in terms of Quillen lifting properties and the ``shift'' endomorphism $[+1]:\sFilth\lra\sFilth$. 
\vskip 4pt

The discussion of limits, Cauchy sequences, equicontinuity suggests the following direction.

\subsubsection{Reformulate compactness and completeness} Develop the theory of compactness and completeness, 
in particular compactness and completeness of function spaces, Arzela-Ascoli theorems, in terms 
of Quillen lifting properties, the shift endomorphism $[+1]:\sFilth\lra\sFilth$, and spaces of maps in $\sFilth$.
\vskip 4pt

\subsubsection*{Acknowledgements} [G] thanks M.Bays for many useful discussions, proofreading and generally 
taking interest in this work. Thanks are also due to  S.Ivanov, S.Kryzhevich, S.Sinchuk, V.Sosnilo, and A.L.Smirnov for a number of discussions,
and particularly to the participants of the seminar of A.L.Smirnov. Importantly, 
thanks are due to D.Rudskii for bringing the Levy-Prohorov metric to our attention. 
We thank D.Grayson for several insightful and encouraging comments.
[G] expresses deep gratitude to
friends for creating an excellent social environment in St.Petersburg, invitations and enabling me to work, 
and the working environment for allowing to freely pursue research interests. 


\end{document}